\documentclass{amsart}
\usepackage{amsfonts}
\usepackage{amsmath}
\usepackage{amssymb}
\usepackage{amsthm}
\usepackage[autostyle]{csquotes}
\usepackage{hyperref}
\usepackage{soul}

\newcommand{\keyw}[1]{\textit{Keywords---} #1}

\title{On first-order arithmetic truth}
\author{Stephen Boyce}
\newtheorem{Theorem}{Theorem}[section]
\newtheorem{Corollary}[Theorem]{Corollary}

\newtheorem{Proposition}[Theorem]{Proposition}
\newtheorem{Lemma}[Theorem]{Lemma}
\theoremstyle{definition}
\newtheorem{Definition}[Theorem]{Definition}

\begin{document}
\begin{abstract}
  The standard interpretation of first-order number
  theory ($PA$), according to the generally accepted view, associates well-defined set-theoretic
  entities with each and every well-formed formula of this system. But this implies that the class of
  $PA$ theorems is semantically defined by a class sign of $PA$ itself, $(\exists x_2) Pf(x_2, x_1)$, in
  the following sense: with $\overline{b}$ the $PA$ numeral for the number $b$, $(\exists x_2) Pf(x_2, \overline{b})$
  is true under the standard interpretation if and only if $b$ is the G\"odel number of a $PA$ theorem. From this
  however it is easily established, by a modification of G\"odel's proof, that the class of $PA$ theorems, and hence the standard interpretation
  of $PA$ itself, is not well defined after all.
\end{abstract}
\maketitle 
\begin{msc}
  03C62,
  03B25,
  03F40,
  03F03
\end{msc} \\
\keyw{Models, First-order number theory, Metamathematics}
\pagenumbering{arabic}
\section{Introduction} \label{section_introduction_truth}
Whilst G\"odel's \cite{godel1931} original incompleteness proof focused on a higher-order system,
with an indication of the application of the results
to each member of a broad class of formal systems, more recent presentations of this proof
often focus on systems like first-order number theory. The wider context of G\"odel's proof has of course changed in various ways since the time of publication.
The widespread acceptance of
Tarski's \cite{tarski1936} theory of truth for formal languages is one such change. Presentations of G\"odel's incompleteness results are now often accompanied
by discussion of semantic issues that were essentially avoided in G\"odel's original presentation. For example, the standard interpretation of first-order arithmetic, or some variation thereof,
may be appealed to to establish that the G\"odel sentence for such a theory asserts of itself that it is not formally provable in the theory in question (\cite{mendelson2015}: 209).

In light of such uses, the question arises: may facts established by
considering the standard interpretation of such systems be used in determining, metatheoretically, whether or not a formula of the system is a theorem? This paper aims to show that if
such facts are used, it follows that that the class of $PA$ theorems, and hence the standard interpretation of $PA$ itself, is not well defined.

As any proof, to avoid an infinite regress, must assume certain premises the prospective reader
is entitled to an indication of the material assumed in what follows. In brief, I take as given
each of the items listed below:
\begin{enumerate}
\item A classical approach to logic.
\item Mendelson's account of first-order number theory $S$ (\cite{mendelson2015}: Chapter 3),
  which, loosely speaking is the theory that results, in the language of (first-order) arithmetic,
  when proper axioms corresponding to Peano's postulates, and the consequences thereof, are
  added to the logical truths for this language.\footnote{The above vague description is not
  meant to be a substitute for a precise definition, including a clarification of the distinction
  between an informal statement of induction for the natural numbers, which quantifies over all
  sets of natural numbers, and the more restricted axiom schemata  ($S9$) of $S$. For
  convenience, at various points below I take $S$ to be the first-order theory $K$ considered by Mendelson in establishing various results.
  At certain points Mendelson's proof (\cite{mendelson2015}: 219) of Proposition 3.43 is followed quite closely.}
\item \label{item_metatheorems} The validity of various metatheorems / derived rules for this system (\cite{mendelson2015}: \S 2.5),
  such as, for example, the deduction theorem (\cite{mendelson2015}: Proposition 2.5). 
\item Standard, Tarskian definitions of various semantic notions, including: interpretation, model, truth under an interpretation (\cite{mendelson2015}: \S 2.2);
\item A specific arithmetization of the syntax of the system, mentioned uniformly throughout the following discussion (\cite{mendelson2015}: \S 3.4);
\item Associated with the above mentioned arithmetization, symbolism for mentioning and using various required arithmetized syntactical relations
  (primarily those mentioned in the following Definition \ref{definitions_notation_symbolism}).
\end{enumerate}

These assumptions are not, on most points, as restrictive as they might seem. Whilst 
I generally follow Mendelson's account (\cite{mendelson2015}) for matters of detail,
familiarity with any classical treatment will do in so far as I have succeeded in
avoiding reliance on anything that is unique to Mendelson's account.

If we adopt Quine's (\cite{quine1970}) view of what he dubbed deviant logics, then the
restriction to a classical approach is similarly inconsequential: the deviant contributions are not
alternative logics but explorations of a different topic, since the newly coined "negations" etc. are not in
fact terms of logic. If, contrariwise, we adopt the Tower of Babel
view of Metamathematics, the restriction is substantial. On this view however the classical focus
facilitates a more in-depth consideration of this perspective.
The would-be encyclopaedic contribution is likely to serve at least one, if not all perspectives poorly.

\section{The class of theorems of first-order number theory}

The demonstration that follows has been divided somewhat artificially into two parts. Firstly, I introduce the
notation and terminology, together with a number of metatheorems and propositions. (For most of the latter items
I simply cite an accepted proof, though on some points I present a somewhat self-contained proof that
utilises the material already at hand.) The second part of the demonstration presents two formal $S$ proofs that make use of
various metatheorems presented in the first part.

\subsection{Notation, terminology and metatheorems}

The following symbolic / notational conventions are used at throughout the following:
\begin{Definition} \label{definitions_notation_symbolism} Notation / symbolism: \newline
  \vspace{-2ex}
  \begin{enumerate}
  \item For any natural number $y$, let $\overline{y}$ name the $S$ numeral for the
    number $y$.
  \item Adapting G\"odel's \cite{godel1931} terminology, an $S$ formula with one (or several) free variable(s), $x$ (or $x, y, \ldots$),
    will sometimes be referred to as a class sign (or relation sign respectively).
  \item If the $S$ relation sign $\mathcal{B}$ contains free occurrences of the $S$ variables $x, y, \ldots$ - which may be indicated by writing $\mathcal{B}$ as
    $\mathcal{B}(x, y, \ldots)$ - then, provided the $S$ terms $a, b, \ldots$ are free for $x, y, \ldots$ respectively, in the usual sense,
    then $\mathcal{B}(a, b, \ldots)$ may be used to name the formula that results when $a, b, \ldots$ respectively are (simultaneously) substituted for all free occurrences of
    $x, y, \ldots$ in $\mathcal{B}(x, y, \ldots)$; where $\mathcal{B}$ does not contain any free occurrences of the variables $x, y, \ldots$,
    then both  $\mathcal{B}(x, y, \ldots)$ and $\mathcal{B}(a, b, \ldots)$ are simply $\mathcal{B}$. (Such a simultaneous substitution may require that
    if there are any occurrences of $x, y, \ldots$ in the terms $a, b, \ldots$ themselves, then the $x, y, \ldots$ to be replaced in $\mathcal{B}(x, y, \ldots)$
    are first replaced by other fresh $S$ variables (\cite{shankar1994}:89).)
  \item Let $\mathfrak{M}$ be the standard interpretation of $S$.
  \item Let "$\vDash_{\mathfrak{M}} \mathcal{B}$"
    assert that the $S$-formula $\mathcal{B}$ is true under the standard interpretation.\footnote{Quine's corners would enable this and the following item to be expressed more precisely,
    but some readers would unfortunately confuse such a use of corners with an alternative usage concerned with the naming of G\"odel numbers.}
  \item Let $T_{s}$ be the class of G\"odel numbers of $S$ theorems.
  \item Let  "$\vdash_{S} \mathcal{B}$" assert $S$ theoremhood, i.e. the existence of an $S$ proof of
    the $S$-formula $\mathcal{B}$.
  \item Conditional proof: Where $\mathcal{B}$ and each member of some finite set $\mathcal{A}_0, \ldots \mathcal{A}_n$
    are all $S$ formulae, $\mathcal{A}_0, \ldots \mathcal{A}_n \vdash_{S} \mathcal{B}$
    asserts that: if proofs of each hypothesis $\mathcal{A}_0, \ldots \mathcal{A}_n$
    exist, then a proof of $\mathcal{B}$ also exists.
  \end{enumerate}
\end{Definition}

With the aid of these symbolic conventions, the content and use of the metatheorems mentioned at Item \ref{item_metatheorems} (in \S \ref{section_introduction_truth})
can be clarified a little further as follows.
\begin{enumerate}
\item Following the metamathematical approach, a formal $S$ proof of an $S$ formula
  $\mathcal{B}$, is a finite sequence
  of $S$ formulae (hereon, always assumed to be well-formed) $\mathcal{A}_0, \ldots \mathcal{A}_n$, such that:
  \begin{enumerate}
  \item Each formula in the sequence is either an $S$ axiom, or follows from one or more
    previous formula(e) in the sequence via an $S$ rule of inference;
  \item $\mathcal{B}$ is $\mathcal{A}_n$ (where $n$ may be zero if $\mathcal{B}$ is an $S$ axiom).
  \end{enumerate}
\item Each of the metatheorems / derived rules used below (in \S \ref{subsection_two_s_proof_lemmas})
  asserts that: if proofs of each member of some finite set of hypotheses $\mathcal{A}_0, \ldots \mathcal{A}_n$
  exist, then a proof of some other formula,  $\mathcal{B}$ also exists;
\item Since the Deduction Theorem is only applied herein to sentences, the simpler formulation framed with respect to the
  propositional calculus may be applied  (\cite{mendelson2015}: Proposition 1.9):  Where $\mathcal{B}$ and each member of some finite set $\mathcal{A}_0, \ldots \mathcal{A}_n$
  are all $S$ formulae, and $\mathcal{A}_0, \ldots \mathcal{A}_n \vdash_{S} \mathcal{B}$ obtains, then
  $\mathcal{A}_0, \ldots \mathcal{A}_{n-1} \vdash_{S} \mathcal{A}_n \Rightarrow \mathcal{B}$ may be inferred.
  For the case when $n$ is zero, we may infer   $\vdash_{S} \mathcal{A}_0 \Rightarrow \mathcal{B}$
  from   $\mathcal{A}_0  \vdash_{S} \mathcal{B}$ (which is referred to below as discharging the hypothesis).
\item With the exception of two rules explained at the following two items, each of the rules are used as discussed at (\cite{mendelson2015}: \S 2.5) -
  with the one minor variation being that, for proof by contradiction, I firstly infer $\tau  \vdash_{S}  (\lnot \mathcal{B}) \Rightarrow  (\mathcal{C}\land \lnot \mathcal{C})$
  from   $\tau , \lnot \mathcal{B} \vdash_{S}  (\mathcal{C}\land \lnot \mathcal{C})$ before inferring  $\tau \vdash_{S} \mathcal{B}$, and allow
  $(\lnot \mathcal{C} \land \mathcal{C})$  in place of  $(\mathcal{C} \land \lnot \mathcal{C})$.
\item At the second step in the proof of Lemma \ref{lemma_one}, to emphasise that a conditional proof is being asserted,
  I cite a metatheorem I have labeled \newline "$\text{Identity: } A \Rightarrow A$",
  to justify the claim of an instance of $\mathcal{B} \vdash_{S} \mathcal{B}$; Mendelson simply proves (\cite{mendelson2015}: Lemma 1.8) that
  "$\vdash \mathcal{B} \Rightarrow \mathcal{B}$" holds for any formula $\mathcal{B}$ and never explicitly introduces the
  derived rule just mentioned.
\item Corollary \ref{corollary_one} introduces a derived rule, proven at Proposition \ref{propostion_theorems_semantically_defined},
  which is possibly used in metatheoretical reasoning about the standard interpretation of $S$ (e.g. \cite{mendelson2015}: 209),
  though this idea may be quite controversial and I strongly suspect Mendelson himself would have objected to it.
\end{enumerate}

For the demonstration of interest, a series of connected propositions concerning the syntax of and standard interpretation of $S$ are required.
\begin{Proposition}
  There exists a (recursive) arithmetized proof relation for $S$ (\cite{mendelson2015}: Proposition 3.28(22)), $\mathrm{Pf}(x, y)$, such,
  for every pair of natural numbers $x$, $y$: $\mathrm{Pf}(x, y)$ holds if and only if (iff)
  $x$ is the G\"odel number of an $S$ proof of the $S$ formula with G\"odel number $y$.
\end{Proposition}
\begin{Proposition} \label{proposition_proof_sign}
  Since every recursive relation may be (numeralwise) expressed - or \emph{syntactically} defined - in $S$ by an
  $S$ relation sign (\cite{mendelson2015}: Corollary 3.31), let $\mathrm{Pf}(x, y)$ be thus syntactically defined in $S$ by the
  $S$ relation sign $Pf(x, y)$, so that, for every pair of natural
  numbers $x$, $y$:
  \begin{equation}
    \mathrm{Pf}(x, y) \rightarrow\ \vdash_{S} Pf(\overline{x}, \overline{y}) \tag{\ref{proposition_proof_sign}.1}
  \end{equation}
  \begin{equation}
    \overline{\mathrm{Pf}}(x, y) \rightarrow\ \vdash_{S} \lnot Pf(\overline{x}, \overline{y}) \tag{\ref{proposition_proof_sign}.2}
  \end{equation}
\end{Proposition}
The following proposition, as mentioned above, is possibly used in metatheoretical reasoning about the standard interpretation of $S$ (e.g. \cite{mendelson2015}: 209),
though I leave it to the reader to form a more definite opinion about such matters:
\begin{Proposition} \label{propostion_proof_semantically_defined}
  The assumption that the standard interpretation, $\mathfrak{M}$, is a model of $S$ implies that
  the the (recursive) arithmetic relation $\mathrm{Pf}(x, y)$ is (truth-functionally) equal to the
  relation associated with the relation sign $Pf(x, y)$ 
  under the standard interpretation in the sense that, for every pair of natural numbers $x$, $y$:
 \begin{equation}
    \mathrm{Pf}(x, y) .\equiv .\ \vDash_{\mathfrak{M}} Pf(\overline{x}, \overline{y}) \tag{\ref{propostion_proof_semantically_defined}.1}
  \end{equation}
\end{Proposition}
  \begin{proof}
    Proofs by contradiction for both the if and only if directions are as follows.
    \begin{description}
    \item[$\mathrm{Pf}(x, y) \rightarrow\ \vDash_{\mathfrak{M}} Pf(\overline{x}, \overline{y})$]
      Assume in order to derive a contradiction that there exists a pair of natural numbers, $x$ and $y$, such that both
      $\mathrm{Pf}(x, y)$ and $\nvDash_{\mathfrak{M}} Pf(\overline{x}, \overline{y})$ hold. From the first hypothesis
      and (\ref{proposition_proof_sign}.1) we may infer \newline
      $\vdash_{S} Pf(\overline{x}, \overline{y})$.
      Since $\mathfrak{M}$ is a model of $S$ this yields $\vDash_{\mathfrak{M}} Pf(\overline{x}, \overline{y})$, a contradiction
      with respect to the second hypothesis.
    \item[$\vDash_{\mathfrak{M}} Pf(\overline{x}, \overline{y}) \rightarrow\ \mathrm{Pf}(x, y)$]
      Assume in order to derive a contradiction that there exists a pair of natural numbers, $x$ and $y$, such that both
      $\vDash_{\mathfrak{M}} Pf(\overline{x}, \overline{y})$ and $\overline{\mathrm{Pf}}(x, y)$ hold. From the second hypothesis
      and (\ref{proposition_proof_sign}.2) we may infer $\vdash_{S} \lnot Pf(\overline{x}, \overline{y})$.
      Since $\mathfrak{M}$ is a model of $S$ this yields $\vDash_{\mathfrak{M}} \lnot Pf(\overline{x}, \overline{y})$, a contradiction
      with respect to the first hypothesis.
    \end{description}
  \end{proof}

  Proposition \ref{propostion_proof_semantically_defined}, which may also be paraphrased as the statement
  that the (recursive) arithmetic relation $\mathrm{Pf}(x, y)$ is, in effect, associated with the relation sign $Pf(x, y)$ 
  under the standard interpretation, is used in the proof of the following proposition.
  \begin{Proposition} \label{propostion_theorems_semantically_defined}
    The assumption that the standard interpretation, $\mathfrak{M}$, is a model of $S$ implies,
    in light of Proposition \ref{propostion_proof_semantically_defined}, that, for any natural number $y$:
    $y$ is a member of the class of G\"odel numbers of $S$ theorems ($T_{s}$) if and only if the sentence ``$(\exists x_2) Pf(x_2, \overline{y})$''
    is true under the standard interpretation (c.f. \cite{swierczkowski2003}: 9):
    \begin{equation}\label{equation_S_Theorems}
      y \in T_{s}\ \text{iff}\ \vDash_{\mathfrak{M}} (\exists x_2) Pf(x_2, \overline{y}) \tag{\ref{propostion_theorems_semantically_defined}.3}
    \end{equation}
  \end{Proposition}
  \begin{proof}
    Proofs of both the if and only if directions are straightforward.
    \begin{description}
    \item[$y \in T_{s}\ \text{if}\ \vDash_{\mathfrak{M}} (\exists x_2) Pf(x_2, \overline{y})$]
      By standard, Tarskian, semantics: if \newline $\vDash_{\mathfrak{M}} (\exists x_2) Pf(x_2, \overline{y})$ holds, there exists a natural number, x,
      such that \newline $\vDash_{\mathfrak{M}} Pf(\overline{x}, \overline{y})$ is also true. Since the relation $\mathrm{Pf}(x, y)$ is the
      relation associated with $Pf(x, y)$ under the standard interpretation, this implies that $x$ is the G\"odel number of an $S$ proof of the formula with G\"odel number $y$.
      From this the truth of $y \in T_{s}$ follows immediately.
    \item[$y \in T_{s}\ \text{only if}\ \vDash_{\mathfrak{M}} (\exists x_2) Pf(x_2, \overline{y})$]
      The hypothesis that $y \in T_{s}$ holds implies that there exists an $S$ proof of the formula with G\"odel number y.
      Let the G\"odel number of this proof be $x$. Hence the relation $\mathrm{Pf}(x, y)$ holds. Since $Pf(x, y)$ \emph{syntactically} defines
      this relation in $S$, we may infer  $\vdash_{S} Pf(\overline{x}, \overline{y})$ and hence, via Existential Rule
      $E4$ (\cite{mendelson2015}: \S 2.5.2), $\vdash_{S} (\exists x_2) Pf(x_2, \overline{y})$ also. We then obtain
      $\vDash_{\mathfrak{M}} (\exists x_2) Pf(x_2, \overline{y})$ in light of the assumption that $\mathfrak{M}$ is a model of $S$.
    \end{description}
  \end{proof}

  Equation \ref{propostion_theorems_semantically_defined}.3, may also be paraphrased as the assertion
  that the arithmetic image of the notion of ``formal provability in $S$'' is \emph{semantically} defined by a formula of $S$ itself
  under the standard interpretation. For ease of reference, I will identify as a Corollary
  a derived inference rule according to which valid inferences may be drawn from Equation \ref{equation_S_Theorems}
  regarding the theoremhood of an $S$ formula.
  \begin{Corollary} \label{corollary_one}
    The (conditional) truth (or falsity) of $(\exists x_2) Pf(x_2, \overline{y})$, under the standard interpretation of $S$, may be used in determining, metatheoretically,
    whether the $S$ formula with G\"odel number $y$ is (or is not respectively) an $S$ theorem.
  \end{Corollary}
  
  A number of standard results also required for the demonstration concern the diagonal function $D(x)$, defined for $S$ as follows.
  \begin{Definition} \label{definitions_diagonal_function} 
    Where $x$ is the G\"odel number of an $S$ formula $\mathcal{B}(v)$ which
    contains free occurrences of the $S$ variable $v$, $D(x)$ is the
    G\"odel number of $S$ formula $\mathcal{B}(\overline{x})$ that results from
    replacing all free occurrences of $v$ in  $\mathcal{B}(v)$ by $\overline{x}$  (c.f. \cite{mendelson2015}: Proposition 3.27(19)).
  \end{Definition}
  The following three results concerning this diagonal function are used:
  \begin{Proposition} \label{propostion_godel_sentence}
    \begin{enumerate}
    \item The diagonal function for $S$ is recursive (c.f. \cite{mendelson2015}: Proposition 3.27(19)).
    \item From Item 1 and \cite{mendelson2015}: Proposition 3.24 we may infer that
      the diagonal function for $S$ is syntactically defined (or representable) in $S$.
    \item Item 2 Implies that the diagonalization lemma (\cite{mendelson2015}: Proposition 3.34) applies to $S$, from which we may infer that
      there exists an $S$-sentence $\mathcal{G}$, the G\"odel sentence for $S$, with G\"odel number $q$ 
      such that:\footnote{Equation \ref{equation_g} is taken from \cite{mendelson2015}: 208, with the modification of substituting the $S$ numeral for $q$, ``$\overline{q}$'',
      in place of ``$\ulcorner \mathcal{G} \urcorner$''.}
      \begin{equation*} \label{equation_g}
        \vdash_S \mathcal{G} \Leftrightarrow (\forall x_2) \lnot Pf(x_2, \overline{q}) \tag{\ref{propostion_godel_sentence}.1}
      \end{equation*}
    \end{enumerate}
  \end{Proposition}

With the above material available we may now consider two $S$ proofs that lead to the main result of interest.

\subsection{Two Lemmas established by $S$ proofs} \label{subsection_two_s_proof_lemmas}

Both of the $S$ proofs to be considered concern the theoremhood of the G\"odel sentence for $S$ defined above.

\begin{Lemma} \label{lemma_one} Corollary \ref{corollary_one}, in conjunction with Equations (\ref{propostion_godel_sentence}.1) and (\ref{propostion_theorems_semantically_defined}.3),
  implies that $(\forall x_2) \lnot Pf(x_2, \overline{q})$ is an $S$ theorem.
\end{Lemma}
To reduce visual clutter in the following proof, the maintained hypothesis \newline
($(\exists x_2) Pf(x_2, \overline{q})$) is not printed to the left of the turnstile
in every line of the deduction prior to its discharge at (8), as, strictly speaking, it ought to be.
\begin{proof}
  \begin{quote}
    \begin{flalign*}
      \quad (1)\ (\exists x_2) Pf(x_2, \overline{q})  && \text{Hypothesis}
    \end{flalign*}
    \begin{flalign*}
      \quad (2)\ (\exists x_2) Pf(x_2, \overline{q}) \vdash_S (\exists x_2) Pf(x_2, \overline{q}) && \text{(1),}\ \text{Identity: } A \Rightarrow A
    \end{flalign*}
    \begin{flalign*}
      \quad (3)\  \vdash_S (\lnot \mathcal{G}) && \text{(2),}\ (\ref{propostion_godel_sentence}.1),\ \text{biconditional elimination}
    \end{flalign*}
    \begin{flalign*}
      \quad (4)\ \vDash_{\mathfrak{M}} (\exists x_2) Pf(x_2, \overline{q})  && (2), \mathfrak{M}\ \text{is a model of $S$}
    \end{flalign*}
    \begin{flalign*}
      \quad (5)\ \vdash_S \mathcal{G} && (4), (\ref{propostion_theorems_semantically_defined}.3) / \text{Corollary}\ \ref{corollary_one}
    \end{flalign*}
    \begin{flalign*}
      \quad (6)\ \vdash_S ((\lnot \mathcal{G})\land \mathcal{G}) && (3), (5)\ \text{conjunction introduction}
    \end{flalign*}
    \begin{flalign*}
      \quad (7)\ (\exists x_2) Pf(x_2, \overline{q}) \vdash_S ((\lnot \mathcal{G})\land \mathcal{G}) && \text{Restatement, (1) - (6)}
    \end{flalign*}
    \begin{flalign*}
      \quad (8)\ \vdash_S [(\exists x_2) Pf(x_2, \overline{q})] \Rightarrow ((\lnot \mathcal{G})\land \mathcal{G}) && (7), \text{deduction theorem}  
    \end{flalign*}
    \begin{flalign*}
      \quad (9)\  \vdash_S \lnot [(\exists x_2) Pf(x_2, \overline{q})] && (8), \text{proof by contradiction} 
    \end{flalign*}
    \begin{flalign*}
      \quad (10)\ \vdash_S (\forall x_2) \lnot Pf(x_2, \overline{q}) && (9), \text{definition of $\exists$, negation elimination} 
    \end{flalign*}
  \end{quote}
\end{proof}

The main result follows easily from the above Lemma.
\begin{Lemma}\label{lemma_two}  Corollary \ref{corollary_one}, in conjunction with Equations (\ref{propostion_godel_sentence}.1), (\ref{propostion_theorems_semantically_defined}.3) and
  Lemma \ref{lemma_one}, implies that the class of $S$ theorems is not well defined.
\end{Lemma}
\begin{proof}
  \begin{quote}
    \begin{flalign*}
      \quad (1)\  \vdash_S (\forall x_2) \lnot Pf(x_2, \overline{q})  && \text{Lemma \ref{lemma_one}}
    \end{flalign*}
    \begin{flalign*}
      \quad (2)\  \vdash_S \mathcal{G} && \text{(1),}\ (\ref{propostion_godel_sentence}.1),\ \text{biconditional elimination}
    \end{flalign*}
    \begin{flalign*}
      \quad (3)\ \vDash_{\mathfrak{M}} (\forall x_2) \lnot Pf(x_2, \overline{q}) && \text{(1),}\ \text{$\mathfrak{M}$ is a model of S}
    \end{flalign*}
    \begin{flalign*}
      \quad (4)\ \nvdash_S \mathcal{G}  && (3), (\ref{propostion_theorems_semantically_defined}.3) / \text{Corollary}\ \ref{corollary_one}
    \end{flalign*}
    \begin{flalign*}
      \quad (5)\ \vDash_{\mathfrak{M}} (\exists x_2) Pf(x_2, \overline{q}) && \text{(2),}\ (\ref{propostion_theorems_semantically_defined}.3) / \text{Corollary}\ \ref{corollary_one}
    \end{flalign*}
  \end{quote}
  $(2)$ and $(4)$, as well as $(3)$ and $(5)$, establish that the class of $S$ theorems is not well defined;
  the semantic metatheory of $S$ is subject to paradox and hence the standard interpretation is also not well defined.
\end{proof}

\section{Discussion}
As the result is at odds with most of what is currently accepted in this area the intrepid reader who has considered it will no doubt
have their own opinion as to its meaning. The scope of the literature to be thus considered is, needless to say, far beyond the aims of this brief note.
More recent discussions of the issues first raised by G\"odel (\cite{godel1931}) often focus on the question of whether various sets are recursively defined,
or the theoretical properties of Turing machines, in a way that tends to sideline considerations of the properties of formal systems (as initially conceived).
Whilst providing additional perspectives to be considered, none of this obviates the need for a more focused consideration, at times, of what can or cannot be established,
by various methods, about formal systems themselves.\footnote{On a
very minor point, I should note that \emph{if} the proof of the theoremhood of $\mathcal{G}$ (Lemma \ref{lemma_one}) is accepted, neither Corollary \ref{corollary_one}
nor Equation (\ref{propostion_theorems_semantically_defined}.3) are required to establish that the that class of $S$ theorems is not well-defined. An alternative, purely syntactical argument may also be
used (see e.g. (\cite{mendelson2015}: Proposition 3.37(a))).}

If the proofs of the two Lemmas are accepted, the skeptical reader might seek to take issue with
Equation \ref{propostion_theorems_semantically_defined}.3 / Corollary \ref{corollary_one}
(and hence also, presumably, reject Proposition \ref{propostion_theorems_semantically_defined}). As best as I can make out, this would amount to
a quite unexpected and substantial revision of material that has been widely accepted for many decades.

The reader whose interest is piqued by the above might also consider a related piece of work, wherein it is shown that a
purely syntactical argument may also be used to establish that the class of $S$ theorems is not well-defined (\cite{boyce2025}). The two
results, taken together, suggests a failure of the devices relied upon for the avoidance of paradox in the metamathematical definition of a formal theory,
rather than a problem confined to the theory of truth for formal languages which is applied. The dovetailing of
the syntactic and semantic approaches on this issue provides strong support for the view that a failure of the basic metamathematical notions of a proof, formal theory etc. is involved.

\section{Conclusion}
G\"odel's (\cite{godel1931}) initial proof of his First Incompleteness Theorem presented, in the informal overview, a \emph{semantic} argument
to support the claim that paradox is avoided:
\begin{quote}
  From the remark that $[R(q);q]$ says about itself that it is not provable, it
  follows at once that $[R(q);q]$ is true, for $[R(q);q]$ is indeed unprovable (being
  undecidable). [\cite{godel1931}: 151]
\end{quote}
As this informal material is put forward only to assist the reader in following the subsequent precise statement of the proof,
it would be unfair to suggest that G\"odel is
committed to defending it. (Which is just as well for G\"odel, as, in light of Tarski's discussion, the prospects for precisely
defining a notion of truth for G\"odel's system $P$, in a way that avoids the liar paradox, seem quite grim (\cite{tarski1936}:\S 5)).
Nevertheless, the tendency has continued since to rely on such semantic arguments to establish that paradox is avoided (e.g. \cite{mendelson2015}: 209).
Whilst other methods of proof may be brought in to play for some systems (e.g. \cite{mendelson2015}: Appendix C), the above demonstration
hopefully establishes, if nothing else, that the intuitions appealed to by these semantic arguments are not as clear and precise as is
commonly thought. If, on the other hand, the exhibited proofs are accepted a far more dramatic revision of the intellectual landscape of
this area is implied.

\end{document}